\documentclass{amsart}

\usepackage{amsmath,amssymb,amsthm,eucal,amscd}

\usepackage[all]{xy}

\theoremstyle{plain}
        \newtheorem{thm}{Theorem}[section]
        \newtheorem{prop}[thm]{Proposition}
        \newtheorem{lemma}[thm]{Lemma}
        \newtheorem{cor}[thm]{Corollary}

\theoremstyle{definition}
        \newtheorem{defn}[thm]{Definition}
        \newtheorem{rem}[thm]{Remark}
 
\DeclareMathOperator{\Aof}{A}
\newcommand{\aof}[1]{\ensuremath{\Aof\left(#1\right)}}

\newcommand{\abof}[2][B]{\ensuremath{\Aof_{#1}\left(#2\right)}}

\newcommand{\parrow}[1]{\quad\downarrow\text{\scriptsize $p_{#1}$}}
\newcommand{\newfib}[1][]
 {\hspace{-6pt}\ensuremath{\underset{#1}{\overset{E_{#1}}{\parrow{#1}}}}
 \hspace{-1pt}}
\newcommand{\aofp}[1][]{\ensuremath{\aof{\newfib[#1]} }}

\newcommand{\aperof}[1]{\ensuremath{\Aof^{\%}\left(#1\right)}}
\newcommand{\sdot}{\mathcal S_{\bullet}}

\DeclareMathOperator{\Qof}{\Sigma^\infty}
\newcommand{\qof}[1]{\ensuremath{\Qof \hspace{-.03in}\left(#1\right)}}

\DeclareMathOperator{\Kof}{K}
\newcommand{\kof}[1]{\ensuremath{\Kof\left(#1\right)}}
\DeclareMathOperator{\Gof}{G}
\newcommand{\gofx}[1][X]{\Gof(#1)}

\newcommand{\atrans}[1]{\tau_{\Aof}\left( #1 \right)}
\newcommand{\atransdec}[1]{\tau_{\Aof}^{dec}\left( #1 \right)}
\newcommand{\ktrans}[1]{\tau_{\Kof} \left( #1 \right)}
\newcommand{\bgtrans}[1]{\tau_{BG} \left(#1\right)}

\DeclareMathOperator{\Retof}{Ret}
\newcommand{\retof}[1][X]{\ensuremath{\Retof^{fd}\left(#1\right)}}
\newcommand{\rof}[1][X]{\ensuremath{\Retof\left(#1\right)}}

\DeclareMathOperator{\Tof}{\mathbb T}
\newcommand{\tof}[1]{\ensuremath{\Tof \left(#1\right)}}
\DeclareMathOperator{\tr}{tr}
\newcommand{\trace}[1][]{\tr_{#1}}
\DeclareMathOperator{\HH}{HH}
\DeclareMathOperator{\THH}{THH}
\newcommand{\thh}[1]{\THH\left( #1 \right)}
\DeclareMathOperator{\MAT}{M}
\newcommand{\mk}[1]{\MAT_k \left( #1 \right)}
\newcommand{\glk}[1]{\widehat {GL}_k(#1)}
\newcommand{\glkn}[1]{{GL}_k(#1)}

\newcommand{\catc}{\mathcal C}
\newcommand{\catd}{\mathcal D}

\DeclareMathOperator{\Nerve}{N}
\newcommand{\nerve}[1]{\Nerve(#1)}
\newcommand{\cynerve}[1]{\Nerve^{cy}(#1)}
\newcommand{\rcynerve}[1]{|\cynerve{#1}|}
\DeclareMathOperator{\Classify}{B}
\newcommand{\classify}[1]{\Classify(#1)}

\DeclareMathOperator{\Map}{map}
\newcommand{\map}{\Map_*}
\newcommand{\mapsp}[1]{\Map \left( #1 \right)}
\DeclareMathOperator{\spectra}{Sp}
\newcommand{\specmap}[1]{\spectra \left( #1 \right)}
\DeclareMathOperator*{\hocolim}{hocolim}
\DeclareMathOperator*{\holim}{holim}
\DeclareMathOperator*{\colim}{colim}
\newcommand{\holimdown}[1]{\underset{#1}{\holim}}
\newcommand{\hocolimdown}[1]{\underset{#1}{\hocolim}}

\newcommand{\Smash}{\wedge}

\DeclareMathOperator{\hocofib}{\ensuremath{hocofib}}

\DeclareMathOperator{\simp}{\ensuremath{simp}}

\newcommand{\diagramit}[1]{\begin{equation*}\xymatrix{#1}\end{equation*}}

\newcommand{\replete}[1][Y]{i_*\retof[#1]}

\newcommand{\pof}[4]{#1 \times #4 \cup_{#2 \times #4} 
  #2 \times #3}

\newcommand{\apof}[4]{\aof{#1 \times #3,\pof{#1}{#2}{#3}{#4}}}

\newcommand{\diff}[1]{Diff(#1)}
\newcommand{\iso}[1]{Iso(#1)}
\newcommand{\piso}[1]{P(#1)}
\newcommand{\spiso}[1]{\mathcal P(#1)}
\newcommand{\swh}[1]{Wh^{diff}(#1)}

\begin{document}
\title[Factoring the Becker-Gottlieb Transfer]
    {Factoring the Becker-Gottlieb Transfer Through the Trace Map}
\author[W. Dorabia{\l}a]
    {Wojciech Dorabia{\l}a}
\address{Department of Mathematics\\
        Penn State Altoona\\
        Altoona, PA 16601-3760}
\email[W. Dorabia{\l}a]{wud2@psu.edu}

\author [M. W. Johnson]
        {Mark W. Johnson}
\email[M. W. Johnson]{mwj3@psu.edu}

\keywords{transfer, trace, algebraic K-theory of spaces, topological Hochschild homology}
\subjclass{Primary: 19D10; Secondary: 18F25, 19Exx, 55R70}

\date{\today}

\begin{abstract}
  In 1998, Becker and Schultz \cite{BecSch} published axioms characterizing the Becker-Gottlieb transfer 
  $\bgtrans{p}: \qof{B_+} \to \qof{E_+}$ for certain types of fibrations $p:E \to B$.
  We verify these axioms
  for the composite of the algebraic $\Kof$-theory transfer $\ktrans{p}:\qof{B_+} \to \aof{E}$ 
  of any perfect fibration $p$
  followed by the evaluation (at the unit) from the free loop space $\Lambda$
  of the B\"okstedt trace map 
  $\trace:\aof{E} \to \qof{\Lambda E_+} \to \qof{E_+}$.  
  As a consequence, for $p$ any compact ANR fibration with finite CW base
  (those considered by Becker-Shultz), $\bgtrans{p} \simeq \trace \ktrans{p}$. 
\end{abstract}

\maketitle

\section{Introduction}

In many situations in mathematics, it is important to understand the collection of automorphisms
of commonly occurring objects.  Thus, much work in the theory of smooth manifolds has focused on understanding $\diff{M}$, the space of self-diffeomorphisms of a compact, smooth manifold $M$.  
Trying to use the tools of homotopy theory leads one to consider the space of (smooth) paths in 
$\diff{M}$ starting at the identity, which is the isotopy space $\iso{M}$.  In 1970, J. Cerf \cite{Cer} introduced the (smooth) pseudoisotopy (or concordance) space $\piso{M}$, consisting of all self-diffeomorphisms of the cylinder $M \times I$ which leave $M \times {0} \cup \partial M \times I$ fixed pointwise.  Notice, by adjointness, an isotopy can also be viewed as a self-diffeomorphism of $M \times I$ which leaves $M \times {0}$ fixed and preserves the second coordinate, so $\iso{M} \subset \piso{M}$.  In fact, whenever $\piso{M}$ is path connected, the more convenient tools of pseudoisotopy are sufficient to study the path components of $\diff{M}$, and Cerf showed that $M$ without boundary, simply connected, and dimension at most 6 implies $\piso{M}$ is path connected.

In order to have more room for homotopies, one next passes to the stable smooth pseudoisotopy space,
$\spiso{M}$, given by the homotopy colimit of the tower of  ``suspension maps"
\[\piso{M} \to \piso{M \times I} \to \piso{M \times I^2} \to \dots  \  .
\]
Igusa's Stability Theorem \cite{Igu} shows the induced stablization map $\piso{M} \to \spiso{M}$ induces an isomorphism for homotopy groups below roughly one third the dimension of $M$.  A delooping 
of $\spiso{M}$ is given by the (smooth) stable h-cobordism space of $M$, which is built from the space of h-cobordisms via ``suspension maps" as well.  A second (or double) delooping 
of $\spiso{M}$, by work of Waldhausen \cite{Wal1}, is given by the (smooth) Whitehead space 
$\swh{M}$.

Using a version of the trace map and topological Hochschild homology,
Waldhausen also showed \cite{Wal3} that there is a natural splitting of infinite loop spaces  
\[
\aof{M} \simeq \qof{M_+} \times \swh{M},
\]
where $\aof{M}$ is his algebraic $\Kof$-theory of spaces, which we refer to as $\Aof$-theory.  
Combined with Igusa's Stability Theorem,
this allows one to recover information about the homotopy groups of $\piso{M}$ (hence indirectly about $\diff{M}$) from information about $\aof{M}$ (and the stable homotopy of the manifold, $\qof{M_+}$).

Another important tool in the study of pseudoisotopy is a geometric transfer map
for pseudoisotopy spaces.  Given $p:N \to M$ a smooth fiber bundle with compact fibers,
this is a map $p^{\ast}:\spiso{M} \to \spiso{N}$, which also passes to the
second delooping $p^{\ast}:\swh{M} \to \swh{N}$.
Burghelea \cite{Bur} 
used these geometric transfer maps for (stable smooth) pseudoisotopy spaces to study the homotopy 
groups of $Diff(M)$.  Under the same condition, 
there are also abstract transfer maps in $\Aof$-theory, 
$\atrans{p}:\aof{M} \to \aof{N}$, which L\"uck \cite{L1}, \cite{L2} verified were 
compatible with the geometric transfers on
$\swh{M}$ via the splitting above, at least for $\pi_0$ and $\pi_1$.  
Thus, one natural question is to understand
the restriction of $\atrans{p}$ to the stable homotopy factor.

Recall that Waldhausen's splitting result has a weaker formulation, saying the composite of 
the trace and the assembly map
\diagramit{
\qof{X_+} \ar[r]^{\alpha} & \aof{X} \ar[r]^-{\trace} & \qof{X_+}  
}
is naturally homotopy equivalent to the identity.   Then, with $\ktrans{p}=\atrans{p} \alpha: \qof{X_+} \to \aof{E}$ denoting the Algebraic $\Kof$-theory transfer, identifying the portion of the $\Aof$-theory transfer visible using stable homotopy, $\trace \atrans{p}|_{\qof{M_+}}$, is equivalent to 
determining $\trace \atrans{p} \alpha=\trace \ktrans{p}$.  
As a consequence of work of Dwyer, Weiss, and Williams 
\cite{DWW}, one sees $\trace \atrans{p}|_{\qof{M_+}}$ is the Becker-Gottlieb(-Dold) transfer for $p$ a smooth fiber bundle.  The primary goal of the current article is to extend this identification for a more general class of fibrations, by exploiting the axiomatic characterization of the Becker-Gottlieb transfer given by Becker and Schultz \cite{BecSch} and stated in our language in Def. \ref{BSaxioms} and Thm \ref{BSthm}.  

Note that in the special case of a disk, $A(D^n) = \qof{S^0} \times \swh{D^n}$ is related to the 
pseudoisotopy (and $h$-cobordism) space of a disc (and the stable homotopy groups of spheres).  
In this (contractible) case, Douglas \cite{Dou} showed $\trace \ktrans{p}$ is the 
Becker-Gottlieb(-Dold) transfer for $p$ any compact ENR-fibration over the disk
as part of a more general result. 
The main result of \cite{Dou} shows that for a fibration $p:E \to X$ with compact fibers making the total space an ENR over the base, the following 
commutes (up to weak homotopy)
\diagramit{
& {\aof{E}} \ar[r]^{\aof{c}} & {\aof{*}} \ar[d]^{\trace} \\
{\qof{X_+}} \ar[r]_-{\bgtrans{p}} \ar[ur]^{\ktrans{p}} & {\qof{E_+}} \ar[r]_{\qof{c_+}} & {\qof{S^0}},
}
where $c:E \to *$ is the collapse map.

The main result here is the following:
\setcounter{thm}{16}
\setcounter{section}{2}

\begin{thm}
	\label{ident}
	If $p:E \to X$ is a compact ANR fibration with finite CW base, then the diagram
\diagramit{
& {\aof{E}} \ar[d]^{\trace} \\
{X} \ar[r]_-{\bgtrans{p}} \ar[ur]^{\ktrans{p}} & {\qof{E_+}} 
}
commutes in the homotopy category of spectra.
\end{thm}

\setcounter{section}{1}

The first author has been investigating, along with Badzioch \cite{BD} as well as both authors with 
Williams \cite{BDW}, applications of this theorem to computations of higher Reidemeister torsion (after Dwyer, Weiss, and Williams \cite{DWW}). 

As mentioned above, the theorem is known for smooth fiber bundles, as a consequence of 
work of Dwyer,Weiss, and Williams \cite{DWW}, where it is shown that for a smooth fiber bundle 
$p:E \to X$, the diagram
\diagramit{
 & {\aof{E}}  \\
{X} \ar[r]_-{\bgtrans{p}} \ar[ur]^{\ktrans{p}} & {\qof{E_+}} \ar[u]_{\alpha}}
commutes up to homotopy.  If one follows each map to $\aof{E}$ with the trace map, 
$\aof{E} \to \qof{E_+}$, Waldhausen's result in \cite{Wal1} that the trace of the assembly map is
homotopic to the identity implies 
\diagramit{
& {\aof{E}} \ar[dr]^{\trace} \\
{X} \ar[ur]^{\ktrans{p}} \ar[r]_-{\bgtrans{p}} & {\qof{E_+}} \ar[u]_{\alpha} \ar[r]_{id} & {\qof{E_+}} 
}
commutes up to homotopy (still in the smooth case) and the lower horizontal is also $\bgtrans{p}$.

Our method of proof for Theorem \ref{main} 
will be to verify that the axioms for characterizing the Becker-Gottlieb transfer given by Becker and Schultz \cite{BecSch} are satisfied by the composite of the trace with the algebraic $\Kof$-theory transfer.  In fact, we will verify that each of these axioms holds (including the stronger version of additivity) for perfect fibrations, that is Hurewicz fibrations whose fibers are retracts up to homotopy of finite CW-cmplexes.  Hence, if it could be shown that these same axioms characterize the Becker-Gottlieb transfer for any larger class of perfect fibrations, 
(consider the preprint of Klein and Williams \cite{Kle}) 
then our result would also extend immediately.

In order to apply the axioms of Becker and Schultz \cite{BecSch}, we must extend to a relative version each of the relevant natural maps, and verify compatibility of these extensions with the natural relative external products along with a relative additivity result.  Since we rely upon maps induced on homotopy cofiber constructions for our relative maps, we cannot simply work in the homotopy category.  We must instead work with weak maps (see subsection \ref{weaknat} below) and study their multiplicative properties carefully.  Fortunately, all of the weak maps of importance to us consist of natural pieces, so a CW pair $(X,Y)$ induces a string of commutative squares (not just homotopy commutative squares, a key detail for us).  As a consequence, there is an induced map on the homotopy cofibers associated to the inclusion at each intermediate stage of the original weak map, and we concatenate these into a relative weak map (see Lemma \ref{pairnat} and Rem. \ref{pairnatrem}).  Of course, this introduces substantial technical complexities, particularly for tracking multiplicative properties, which we address in the earlier sections.

Section 2 is devoted to background material familiar to the experts in this area, along with our definitions for the (relative) $\Aof$-theory transfer (following Williams \cite{Wil}), assembly (following Weiss and Williams \cite{WeiWill}) and trace maps (following B\"okstedt \cite{BHM}).  The section ends with a statement of the axioms and main result of Becker and Schultz (Def. \ref{BSaxioms} and Thm \ref{BSthm}) in the current language, as well as our primary result (Thm \ref{main}).  Section 3 is focused on the rather involved details of the compatibility of the earlier constructions with external pairings, and includes the verification of the multiplicativity (Prop. \ref{prodtransfer})
of our candidate transfer.  The short section 4 is then devoted to verifying the normalization (Prop. \ref{normalization}) and naturality (Prop. \ref{naturality}) properties for our candidate transfer, relying mainly on details presented earlier.  

Section 5 is dedicated to verifying the strong additivity axiom, which is more complicated.  The key
underlying result is the Transfer Additivity Theorem of the first author \cite{Do} and \cite{BD}.  Unfortunately, that result does not include a naturality statement for the relevant homotopy, hence it need not extend to an induced map on homotopy cofibers.  
This requires us to work with a less well-known model for the algebraic $\Kof$-theory transfer, following Dwyer, Weiss, and Williams \cite{DWW}, which is defined in terms of parametrized $\Aof$-theory, and closer to that used by Douglas \cite{Dou}.

\section{Notation and Background}

This section will introduce material already familiar to experts in this area, along with some notation.  Included are various definitions for retractive spaces, $\Aof$-theory, the bivariant $\Aof$-theory introduced by Williams, and a relative version of his view of the $\Aof$-theory transfer, which was inspired by unpublished work of Waldhausen.
This is followed by a careful discussion of weak natural transformations, along with the relative assembly and trace maps.  The section ends with a restatement of the axioms of Becker and Schultz, along with their main theorem, and the statement of our main result.

\subsection{Classes of Fibrations}
\label{allfibs}

First, there are several different classes of (Hurewicz) fibrations (usually denoted $p:E \to X$)
considered in our references.  We
would like to work with the class of perfect fibrations following Williams \cite{Wil}, that is, fibrations whose fibers are retracts
up to homotopy of finite CW complexes.  This class will contain all of the others, so our results will apply to each of the other classes as well.  We will also assume the base space $X$ has the homotopy type of a CW complex.  Among other things, this means $\qof{X_+}$ will have the homotopy type of a CW spectrum, so the Hurewicz Theorem will imply any weak homotopy equivalence between spectra of this type will actually be a homotopy equivalence.  Thus, we avoid the necessity of working with weak homotopy equivalences in many cases.

Given a CW inclusion $i:Y \to X$, a perfect fibration of pairs
$ p=(p_X,p_Y):(E_X,E_Y) \to (X,Y),$ consists of a
perfect fibration $p_X$ and its restriction to the subspace $Y$ 
	(which is then also a perfect fibration).

Becker and Schultz \cite{BecSch} work with compact ANR fibrations when establishing their axioms characterizing the Becker-Gottlieb transfer 
$\bgtrans{p}:X \to \qof{E_+}$.  That is, they consider fiber bundles $p:E \to X$ with $X$ a finite 
CW complex and fibers finite dimensional compact ANRs (Absolute Neighborhood Retracts).  As they emphasize, an important consequence is that $E$ then also has the homotopy type of a finite 
CW complex.  

A different class of fibrations, considered by Douglas \cite{Dou}, are the compact ENR fibrations, that is fibrations with compact fibers that make the total space an ENR (Euclidean Neighborhood Retract) over the base.

\subsection{Retractive Spaces}
\label{retrdef}

A retractive space over $X$ consists of a space $Z$, together with maps 
$i:X \to Z$ and $r:Z \to X$ such that $r \circ i$ is the identity on $X$.  Clearly this becomes a 
category $\rof$ if as morphisms from $(Z,i,r)$ to $(Z',i',r')$ one takes
all (continuous) maps $f:Z \to Z'$ with $f \circ i=i'$ and $r=r' \circ f$.  It is common to leave
the choice of $i$ and $r$ out of the notation whenever possible without confusion.

Given $M \in \rof$ and $N \in \rof[W]$, their external coproduct is the element of $\rof[X \times W]$
defined by the pushout diagram (with both maps coming from the structural inclusions)
\diagramit{
{X \times W} \ar[d] \ar[r] & {X \times N} \ar[d] \\
{M \times W} \ar[r] & {M \times W \cup_{X \times W} X \times M}.
}

We will use $M \hspace{.025in} _X\hspace{-.05in}\vee_{W} N$ as a short form of $M \times W \cup_{X \times W} X \times M$.

Similarly, the external (smash) product is the element of $\rof[X \times W]$ given by the pushout
diagram
\diagramit{
{M \hspace{.025in} _X\hspace{-.05in}\vee_{W} N} \ar[d] \ar[r] & {M \times N} \ar[d] \\
{X \times W} \ar[r] & {M \times N \cup_{M \times W \cup_{X \times W} X \times M} X \times W}.
}

We will use $M\hspace{.025in} _X\hspace{-.05in} \wedge_W N$ as a short form of 
$M \times N \cup_{M \times W \cup_{X \times W} X \times M} X \times W$.

A retractive space $Y$ over $X$ is homotopy finite if there is a morphism of retractive spaces $W \to Y$ which is an underlying homotopy equivalence, such that the structure map $X \to W$ is the inclusion of a finite relative CW complex.  One calls $Z$ (homotopy) finitely dominated
over $X$ if there exists a composition in the category of retractive spaces $Z' \to Y \to Z$ whose underlying composite is a homotopy equivalence of spaces, such that $Y$ is homotopy finite.
We will use $\retof$ to indicate
the full subcategory of retractive spaces over $X$ with objects the finitely dominated retractive spaces.
This is a Waldhausen category (that is, a category with cofibrations and weak equivalences in Waldhausen's terminology \cite{Rut}) where the cofibrations are the closed embeddings satisfying the homotopy extension property (the closed cofibrations of \cite{Str}) and weak equivalences are the homotopy equivalences (see \cite[II.6.1]{DWW} for details).  Thus, we can apply
Waldhausen's $\Kof$-theory functor, which yields the spectrum $\aof{X}$ in this case.

\subsection{Bivariant $\Aof$-theory}
\label{bivardef}

More generally, given $p:E \to X$, one calls
a retractive space over $E$ fiberwise finitely dominated over $p$ if each homotopy fiber of $p \circ r$
at $x \in X$ is finitely dominated over the fiber of $p$ at $x$, $F_x$.  We will use $\retof[p]$ to indicate the full subcategory of retractive spaces over $E$ with objects that are 
fiberwise finitely dominated over $p$.  This is a sub-Waldhausen category of $\retof[E]$, so we can again apply Waldhausen's $\Kof$-theory functor, which yields the spectrum $\aof{p}$ in this case, the Bivariant $\Aof$-theory of Williams \cite{Wil}.  One particularly important thing to notice here is that 
$p:E \to X$ a perfect fibration implies $E \coprod E \to E \to X$ is fiberwise finitely dominated over $p$.  This yields a point $\chi(p)$ in $\aof{p}$ which Williams \cite{Wil} uses to give a particularly useful description of the transfer map in the algebraic $\Kof$-theory of spaces, also exploited by the first author alone 
\cite{Do}, with Badzioch \cite{BD}, and with the second author \cite{DJ}.  We will later refer to this point in the $\Aof$-theory spectrum as the homotopy parametrized Euler characteristic of the fibration.

A CW-inclusion $i:Y \to X$ induces an exact functor $i_*:\retof[Y] \to
\retof$ and we use $\replete$ to denote its replete image, that is the full
subcategory of $\retof$ consisting of objects isomorphic to some
$i_*(W)$.  If we give $\replete$ the Waldhausen structure it inherits
from $\retof$, then $i_*$ even becomes an equivalence of Waldhausen
categories $\retof[Y] \to \replete$, since $i$ is a CW-inclusion by
assumption.   Also, recall that a ``quotient'' of
retractive spaces over $X$ is defined by the pushout
\diagramit{
B_i \ar[r] \ar[d] & B_j \ar[d] \\
X \ar[r] & {B_j/B_i}.
}

\subsection{Constructing the (Relative) $\Aof$-theory transfer}
\label{makerelatrans}

Given a perfect fibration $p:E \to X$, it follows from \cite[Lemma 3.1]{And} 
that finite domination is preserved under pullback over the fibration.  
As a consequence, one has a pullback
functor $\retof \to \retof[E]$ which is exact,
as is readily verified using the following Lemma due to L\"uck.  An earlier formulation of this result is available at the end of \cite{Str}.
 
\begin{lemma}[Lemma 1.26 of L\"uck \cite{Lu}]
  \label{wolfgang}
  Suppose 
\diagramit{
{A} \ar[r]^{f} \ar[d]_{j} & {Y} \ar[d]^{J} \\
{X} \ar[r]_{F} & {Z}
}
  is a pushout diagram while $j$ is a (closed) cofibration.  Given a fibration 
  $p:E \to Z$, the pullback construction 
  yields a pushout diagram
\diagramit{
{f^*J^*E} \ar[r]^{\overline f} \ar[d]_{\overline{j}} & {J^*E} \ar[d]^{\overline J} \\
{F^*E} \ar[r]_{\overline F} & {E}
}
  with $\overline{j}$ also a (closed) cofibration.
\end{lemma}

As an exact functor, this pullback construction induces a map on $\Aof$-theory, known as the
$\Aof$-theory transfer.  We will generally use notation like $p_X^*$ (or even $p^*$ when no confusion will result) for such transfers.  Given a perfect fibration of pairs (see subsection \ref{allfibs}), there
is a commutative diagram of Waldhausen categories and exact functors
\diagramit{
\retof[Y] \ar[r]^{i_*} \ar[d]_{p_Y^*} & \retof \ar[d]^{p_X^*} \\
\retof[E_Y] \ar[r]_{\bar i_*} & \retof[E_X] ,
}
since $p_Y$ is simply the restriction of $p_X$ (and the singleton is assumed to be unique).  
As a consequence, the induced diagram
\diagramit{
\aof{Y} \ar[r]^{i_*} \ar[d]_{p_Y^*} & \aof{X} \ar[d]^{p_X^*} \\
\aof{E_Y} \ar[r]_{\bar i_*} & \aof{E_X}
}
commutes as well.  

Let $\hocofib$ indicate the Bousfield-Kan \cite{BK} model for the homotopy cofiber of spectra.

\begin{defn}
	\label{makepairs}
	Given a functor $G$ from spaces to spectra, define its {\it extension to CW pairs} by
\[G(X,Y)=\hocofib(G(i):G(Y) \to G(X))
\] 	 
	for $i:Y \to X$ a CW inclusion.
\end{defn}

Then we have our relative $\Aof$-theory transfer by functoriality of $\hocofib$.

\begin{lemma}
	\label{transdef}
	For any perfect fibration of pairs $(p_X,p_Y):(E_X,E_Y) \to (X,Y)$, 
	there is an induced relative $\Aof$-theory transfer
	$\atrans{p_X,p_Y}:\aof{X,Y} \to \aof{E_X,E_Y}$.
	Furthermore, given a fiberwise homotopy equivalence $f:E_W \to E_X$ making
	the diagram	
\diagramit{
{(E_W,E_Z)} \ar[r]^{f} \ar[d]_{p_2} & {(E_X,E_Y)} \ar[d]^{p_1} \\
{(W,Z)} \ar[r]_{g} & {(X,Y)}
}
commute, the naturality diagram	
\diagramit{
{\aof{W,Z}} \ar[r]^{\aof{g}} \ar[d]_{\atrans{p_2}}
         & {\aof{X,Y}} \ar[d]^{\atrans{p_1}} \\
{\aof{E_W,E_Z}} \ar[r]_{\aof{f}} & {\aof{E_X,E_Y}}
}
commutes up to homotopy.	
\end{lemma}

\begin{proof}
	The existence follows from functoriality of $\hocofib$ and the commutative diagram
	above Def. \ref{makepairs}.  Now we decompose the naturality result into three special cases,
	by first factoring $g$ as a cofibration $k$ followed by a homotopy equivalence $q$.
	Then define $p_3$ to be the pullback of $p_2$ over $q$, and $p_4$ the pullback
	of $p_3$ over $k$.  This splits the original diagram into three pieces,
\diagramit{
{(E_W,E_Z)} \ar[r]^{h} \ar[dr]_{p_2} & {(P_W,P_Z)} \ar[d]^{p_4} \ar[r]^{\bar k} 
	& {(P_U,P_V)} \ar[d]^{p_3} \ar[r]^{\bar q} & {(E_X,E_Y)} \ar[d]^{p_1} \\
& {(W,Z)} \ar[r]_{k} & {(U,V)} \ar[r]_{q} & {(X,Y)}
}	
	with both $p_3$ and $p_4$ also perfect fibrations of pairs, while $\bar k$, $\bar q$ and
	consequently the induced map $h$ are fiberwise homotopy equivalences.

	First, we consider the case where $p_2$ is the pullback of $p_1$ and $g$ is a 
	closed cofibration.
	Then the pullback construction will yield a cubical diagram of Waldhausen
	categories and exact functors.  Two faces of this cube will commute
	on the nose as above, two more will commute on the nose if we choose models for 
	pushouts (but otherwise only up to natural isomorphism compatible with restriction
	to subspaces),
	and the remaining two commute up to natural isomorphism (still compatible with 
	restriction to subspaces) using \ref{wolfgang}.
	As a consequence, the induced square of spectra, after applying the $\Kof$-theory functor 
	to the cube and $\hocofib$ along the relevant edges, will commute up to homotopy. 	 
	
	Now, consider the case where $p_2$ is the pullback of $p_1$, and $g$ is a homotopy
	equivalence.  Recall $f$ is a fiberwise homotopy equivalence so
	both $1 \to g^*g_*$ and
	$f_*f^* \to 1$ are natural homotopy equivalences on the full categories of
	retractive spaces.  As a consequence we have
\[f_*p_2^* \to f_* p_2^* g^*g_*=f_*f^*p_1^*g_* \to p_1^*g_*
\]
	a natural homotopy equivalence between the relevant exact functors which is also
	compatible with restriction to subspaces.  Hence, as above,
	the induced diagram on relative $\Aof$-theory commutes up to homotopy.
	
	Finally, we consider the case where $g$ is the identity.
	In this case, there is a 
	triangular prizm of Waldhausen categories and exact functors which leads to the 
	triangular diagram of spectra
\diagramit{
\aof{E_W,E_Z} \ar[rr]^{\aof{f}} && \aof{E_X,E_Y} \\
& \aof{W,Z} \ar[ul]^{p_2^*} \ar[ur]_{p_1^*} .
}	
	Since $f$ is a fiberwise homotopy
	equivalence, it again follows that $f_*f^* \to 1$ is a natural homotopy equivalence as 
	functors on the full categories of retractive spaces.  As $p_1 \circ f=p_2$, this leads to a
	natural homotopy equivalence $f_* \circ p_2^* = f_*f^*p_1^* \to p_1^*$,
	also compatible with restriction to the subspace $Z$ as desired.  Once again, this leads to 
	homotopy commutativity of the triangle of spectra, and combined with the previous cases
	completes the proof.
\end{proof}

\subsection{Weak Natural Transformations}
\label{weaknat}

We will be considering various constructions commonly viewed as equivalent in the homotopy category of spectra, also called the stable homotopy category, but we must do so with great care in order to allow us to be precise about induced maps on homotopy cofibers.  As a consequence, we cannot simply work in the homotopy category of spectra, where many of these constructions are usually defined.  We would prefer to use a symmetric monoidal model for the stable homotopy category, such as the S-modules of Elmendorf et al \cite{EKMM} or the Symmetric Spectra of Hovey et al \cite{HSS}, but will generally not need to be involved with any particular choice.  

Following common usage, we have several definitions.

\begin{defn}
	\label{weakmaps}
\begin{itemize}
\item[a)]	Given two spectra (or spaces) $Y$ and $X$, a {\it weak map} $Y \to X$ will consist of a string
		of morphisms of spectra (or spaces)
\diagramit{
Y=X_0 \ar[r]^{f_1} & X_1 & \ar[l]_{f_2} X_2 \ar[r]^-{f_3} & \dots X_n=X, 
} 
	where each (``wrong way") map $f_{2k}$ is a weak equivalence.
\item[b)]	The {\it weak composition} of two weak maps will be defined by concatenation (possibly
	with an identity map added to the second string to maintain our convention that even maps
	are weak equivalences), even though weak maps do not form a category due to 
	cardinality issues.  Note that weak maps will determine morphisms in the homotopy
	category of spectra, and their weak compositions will descend to composition
	in the homotopy category.
\item[c)]	Two weak maps will be called {\it weakly equivalent} precisely when they induce
	the same morphism in the homotopy category of spectra.
\item[d)]	Given two spectrum-valued functors $F$ and $G$, a {\it weak natural transformation}
	$F \to G$ will similarly indicate a string of natural transformations
\diagramit{
F=G_0 \ar[r]^{\eta_1} & G_1 & \ar[l]_{\eta_2} G_2 \ar[r]^-{\eta_3} & G_3 \dots G_n=G, 
} 
	where each $\eta_{2k}$ is a natural weak equivalence. 
\item[e)]	As above, one defines the {\it weak composition} of two weak natural transformations
	by concatenation, which descends to composition of natural transformations into
	the homotopy category of spectra.	
\item[f)]	Once again, two weak natural transformations will be called {\it weakly equivalent}
	precisely when they induce the same morphism in the homotopy category of spectra.
\end{itemize}
\end{defn}
	
\begin{rem}
	\label{weakenough}
It follows from the definition that verifying two weak natural transformations $\eta_n, \dots ,\eta_1$ and $\nu_n, \dots ,\nu_1$ satisfy the condition that each $\eta_i$ is naturally (weakly) homotopy equivalent to $\nu_i$ will suffice to prove the weak natural transformations are weakly equivalent.
\end{rem}

When we speak of a relative weak natural transformation below for functors from CW pairs to spectra, the naturality will be with respect to morphisms of CW pairs, which will be defined as commutative squares in topological spaces
\diagramit{
Y \ar[d]_i \ar[r] & Y' \ar[d]^{i'} \\
X \ar[r] & X'
}
with both $i$ and $i'$ CW inclusions.

\begin{lemma}
	\label{pairnat}
	Given a weak natural transformation of (space or) spectrum valued functors $\eta:F \to G$, there is
	an induced relative weak natural transformation of their extensions to CW pairs,
	$\Psi(\eta):F \to G$.
\end{lemma}

\begin{proof}
	Given a string of natural transformations
\diagramit{
F=G_0 \ar[r]^{\eta_1} & G_1 & \ar[l]_{\eta_2} G_2 \ar[r]^-{\eta_3} & \dots G_n=G, 
} 
	we must produce a similar string for the extensions to CW pairs.  Consider the diagram
\diagramit{
F(Y)=G_0(Y) \ar[d] \ar[r]^-{\eta_1} & G_1(Y)\ar[d] & \ar[d] \ar[l]_{\eta_2} G_2(Y) \ar[d] \ar[r]^-{\eta_3} &
	\dots G_n(Y)=G(Y) \ar[d] \\
F(X)=G_0(X) \ar[d] \ar[r]^-{\eta_1} & G_1(X)\ar[d] & \ar[d] \ar[l]_{\eta_2} G_2(X) \ar[d] \ar[r]^-{\eta_3} &  
	\dots G_n(X)=G(X) \ar[d] \\
F(X,Y)=G_0(X,Y) \ar@{.>}[r]^-{\eta_1} & G_1(X,Y) & \ar@{.>}[l]_{\eta_2} G_2(X,Y) 
	\ar@{.>}[r]^-{\eta_3} &  \dots G_n(X,Y)=G(X,Y) 
}	
	so each column is a homotopy cofiber sequence.
	Since the upper squares each commute by the naturality assumption on each $\eta_i$,
	the dotted arrows exist by the functoriality of this model for the homotopy cofiber.  
	Furthermore, since $\eta_{2k}$
	is assumed to be a natural weak equivalence, we see the same remains true 
	on the bottom row as well, which completes the 
	definition of $\Psi(\eta)$.
\end{proof}

\begin{rem}
	\label{pairnatrem}
	The key in the previous lemma is that weak natural transformations consist of pieces 
	which are natural on the nose, rather than up to homotopy.  Thus, when forming each
	homotopy cofiber vertically, there is no dependence on a choice of homotopy
	when forming the induced map.  
	As a consequence, there is an induced relative weak map which also does not depend on any
	choice of homotopy.
\end{rem}

\subsection{Constructing the (Relative) Assembly Map}
\label{makerelass}

Our first example of a weak natural transformation of spectrum valued functors of spaces
will be the assembly map $\qof{X_+} \to \aof{X}$.
Here we consider the composite
\diagramit{
\qof{X_+} \ar[r]^-{\simeq} & X_+ \Smash S^0 \ar[r]^-{1 \Smash \mu} & X_+ \Smash \aof{*} 
& \ar[l]_-{\simeq} \aperof{X} \ar[r]^-{\alpha}& \aof{X}
}
where $\aperof{X}$ is defined as a certain homotopy colimit, following Weiss and Williams \cite{WeiWill}, and 
$\alpha$ represents their assembly map, while $\mu:S^0 \to \aof{*}$ is the unit map of this ring spectrum.  The naturality claim follows from the fact that for a CW inclusion $i:Y \to X$, each square in the diagram
\diagramit{
\qof{Y_+} \ar[r]^-{\simeq} \ar[d] & Y_+ \Smash S^0 \ar[r]^-{1 \Smash \mu} \ar[d] & Y_+ \Smash \aof{*} 
\ar[d] & \ar[l]_-{\simeq} \aperof{Y} \ar[r]^-{\alpha}\ar[d] & \aof{Y} \ar[d] \\
\qof{X_+} \ar[r]^-{\simeq} & X_+ \Smash S^0 \ar[r]^-{1 \Smash \mu} & X_+ \Smash \aof{*} 
& \ar[l]_-{\simeq} \aperof{X} \ar[r]^-{\alpha}& \aof{X} 
}
commutes.

This allows us to introduce a relative assembly map.  

\begin{cor}
	\label{relassem}
	There is a relative weak natural transformation of spectrum valued functors of CW pairs
\[\Psi(\alpha):\qof{X,Y} \to \aof{X,Y}
\]
	induced by the assembly map described above, via Lemma \ref{pairnat}.
	(See Remark \ref{pairnatrem}.)
\end{cor}

\subsection{FSP Version of $\Aof$-theory}
\label{prodtrick}

In constructing out trace map, we will use another model for the $\Aof$-theory of a (connected) space based on FSPs (functors with smash products, see \cite{BHM} or \cite{Bok}), which comes from Waldhausen's ``ring up to homotopy" approach.  Recall that Kan originally constructed \cite[\S 9]{Kan}, a functor from pointed connnected simplicial sets to 
(free) simplicial groups, the loop group, which comes equipped with a universal bundle.
We will use Waldhausen's modification of the construction \cite{WalK} which, among other things, makes clear that it preserves products.
Hence, applying geometric realization to (the underlying simplicial set of)
the loop group of the singular set of a pointed connected space $X$ defines a functor $\gofx$ 
to topological groups (equipped with a universal $\gofx$ bundle).  
Thus, for a connected space $X$, one chooses a basepoint and then
forms the $\Kof$-theory of the FSP associated to this topological monoid, 
$L_{\gofx}(Z) = \gofx_+ \Smash Z$.   For us, $\Kof$-theory of an FSP is defined as
\[\kof{F}={\mathbb Z} \times ( \hocolim_{k} \classify {\glk{F}})^+.
\]
Notice this choice of basepoint is not relevant in determining the homotopy type of the resulting $\Kof$-theory space, since the loop groups for different choices of basepoint will remain isomorphic (see \cite{WalK}), although the specific isomorphism depends on a choice of path between the two choices of basepoint.  Thus,  we will generally only have to be careful to make compatible choices of basepoints when working with maps.  Of course, we are primarily interested in naturality in spaces to apply Lemma 
\ref{pairnat}, so only for the inclusion of a CW pair, where basepoints can be chosen within the subspace, so this will cause no difficulty.

Notice that if $X=\coprod_\alpha X_\alpha$ is a decomposition into components, then the retractive spaces model has an isomorphism 
$\aof{X} \approx \prod_\alpha \aof{X_\alpha}$ (see e.g. \cite[proof of 2.1.7]{Rut}).  This suggests that one extend the FSP model to arbitrary spaces by simply forming the product of the values on the components indexed on the set of components, which will then remain equivalent to the retractive spaces model even if $X$ is not connected. 

The extended FSP model (for not necessarily connected spaces) will also be natural (with appropriate choice of a basepoint for each component), since one can extend the functoriality of the FSP model over this product construction.  In particular, for $Y \to X$ a CW inclusion, even if $Y$ has more components that $X$, one has a specific map
$\aof{Y} \to \aof{X}$ which is compatible (up to homotopy) with the comparable map on the retractive spaces model.  

In what follows, we will require some notation.  Given an FSP $F$, recall the matrix $\mk{F}$ FSP defined by 
\[\mk{F}(Z) = \map([k],[k] \Smash F(Z))
\]
and let $\mathcal H^n_k (\gofx)$ denote the simplicial monoid of pointed
$\gofx$-equivariant self homotopy equivalences of $[k] \Smash S^n \Smash \gofx_+$.

\begin{prop}
	\label{FSPmodel}
	There is a weak natural transformation of spectra, where each component is a weak equivalence,
	$\aof{X} \to \kof{L_{\gofx}}$.  
\end{prop}

\begin{proof}
There is a natural weak equivalence \cite[Prop 2.1.4]{Rut}
$\aof{X} \to \Omega|hR_{hf}(*,\gofx)|$ and a natural chain of homotopy equivalences 
\cite[Prop. 2.2.1 and below]{Rut} between $\Omega|hR_{hf}(*,\gofx)|$ and 
$\mathbb Z \times \lim_{n,k} \classify{\mathcal H^n_k (\gofx)}^+$.   

Following Dundas, Goodwillie, and McCarthy \cite[\S 2.3.4]{DGM}, 
we first construct a natural weak equivalence
$(\glk{L_{\gofx}})^q \stackrel{\sim}{\to} \hocolim_{x \in I} ({\mathcal H^x_k (\gofx)})^q$ which induces an 
entrywise weak equivalence of simplicial spaces, and thereby a weak equivalence on realizations
$\classify{\glk{L_{\gofx}}} \stackrel{\sim}{\to}  \hocolim_{x \in I} \classify{\mathcal H^n_k (\gofx)}$.
The map is constructed by first commuting (unpointed) $\hocolim_I$ with the product
\diagramit{
(\hocolim_{x \in I} \Omega^x \mk{L_{\gofx}}(S^x))^q \ar[d] \\ 
\hocolim_{(x_1,\dots,x_q) \in I^q} (\Omega^{x_1} \mk{L_{\gofx}}(S^{x_1}) \times \dots \times 
	\Omega^{x_q} \mk{L_{\gofx}}(S^{x_q}))
}
then including each factor in a copy of the larger stable matrix value via appropriate suspensions
\[i_j:\Omega^{x_j} \mk{L_{\gofx}}(S^{x_j}) \to 
	\Omega^{(x_1,\dots,x_q)} \mk{L_{\gofx}}(S^{(x_1,\dots,x_q})
\]
followed by the map
\[\hocolim_{I^q} \Omega^{(x_1,\dots,x_q)} \mk{L_{\gofx}}(S^{(x_1,\dots,x_q}) \to
	\hocolim_I \Omega^x \mk{L_{\gofx}}(S^x)
\]
induced by the natural map $I^q \to I$.  Restricting to units up to homotopy then induces the
requisite natural weak equivalence.  
Now we apply $\mathbb Z \times (\hocolim_k ?)^+$ to get a map
$\kof{L_{\gofx}} \to \mathbb Z \times (\hocolim_k \hocolim_{x \in I} \classify{\mathcal H^x_k (\gofx)})^+$
which we extend to $Z \times \lim_{n,k} \classify{\mathcal H^n_k (\gofx)}^+$
via a natural weak map where each component is
a weak equivalence,  
$\hocolim_k \hocolim_{x \in I} \classify{\mathcal H^x_k (\gofx)} \to
\colim_{k,n} \classify{\mathcal H^n_k (\gofx)}$
\end{proof}

\subsection{Constructing the (Relative) Trace Map}
\label{makereltrace}

The motivation for our trace map comes from the Dennis trace map in the Algebraic $\Kof$-theory of rings.  This is induced by a composite natural transformation
\diagramit{
{\classify{\glkn{R}}} \ar[r]^-{i} & \rcynerve{\glkn{R}} \ar[r]^{m} & \HH(M_n(R)) \ar[r]^-{\HH(tr)} & \HH(R).
}
Here $i$ is a canonical section 
$i:\classify{\glkn{R}} \to \rcynerve{\glkn{R}}$ of the natural projection (forgetting the zero simplex)
$\pi:\rcynerve{\glkn{R}} \to \classify{\glkn{R}}$, 
with $i$ given on simplices by $(h_0, \dots,h_k) \mapsto 
((h_0 \dots h_k)^{-1},h_0, \dots,h_k)$, while $m$ effectively just rewrites tuples as strings of 
tensor products, and
$tr:M_n(R) \to R$ is the usual trace of a matrix.  This yields a map
$\coprod_{k \geq 0} \classify{\glkn{R}} \to \HH(R)$, which (after adjusting $\pi_0$), extends over the group completion to yield the Dennis trace map $\kof{R} \to \HH(R)$.
Unfortunately, the analog of $i$ is complicated in the case of rings up to homotopy by the fact that one must be careful when trying to invert elements which are only invertible up to homotopy.

The B\"okstedt trace map \cite{Bok}, also published in B\"okstedt, Hsiang, and Madsen
\cite{BHM}, defines a
weak natural transformation of spectra (for connected spaces) from 
$\kof{L_{\gofx}}$ to $t\thh{L_{\gofx}}$, which is weakly equivalent to
the stable free loop space $\qof{\Lambda X_+}$, and we will describe it in a moment.  However,
choosing the basepoint $1 \in S^1$ (modeled by the units of $\mathbb C$) 
gives a natural map, evaluation at $1$,
$eval_1:\qof{\Lambda X_+} \to \qof{X_+}$ and the (weak) composite, referred to here as the evaluation 
of the B\"okstedt trace, is another weak natural transformation
 $\kof{L_{\gofx}} \to \qof{X_+}$ (for connected spaces).  For non-connected spaces, we will proceed with a product construction as discussed in subsection \ref{prodtrick}, keeping in mind that the natural map from a finite coproduct to the product of the same spectra is a stable homotopy equivalence.  See also Madsen's discussion \cite{Madsen}.

For the sake of brevity, we will follow \cite{BHM} and focus on the zero spaces in constructing the
B\"okstedt trace map, citing compatibility with the relevant $\Gamma$-space structure to promote to the level of spectra.  Recall that the space $\thh{F} = |\THH_\bullet (F)|$, where the space of $p$-simplices of this simplicial space is
\[\hocolim_{I^{p+1}} \map(S^{i_0} \Smash \dots \Smash S^{i_p},
	F(S^{i_0}) \Smash \dots \Smash F(S^{i_p}))
\]
with face maps associated to concatenation $I^{p+1} \to I^p$ 
(preceded by cyclic permuation in the last case).

Recall that for an FSP $L$, its $k$th general linear monoid, $\glk{L}$, is defined as the group-like
topological monoid 
\[\hocolim_I \map(S^i,\mk{L}(S^i))^{\ast}
\]
where the superscript $\ast$ indicates taking the subspace of homotopy units.

As usual,
for a topological monoid $H$, $\classify {H}$ will refer to the realization of the 
simplicial space $\nerve{H}$
whose space of $p$-simplices is given by $p$-fold product of copies of $H$, with face maps given by multiplication of neighboring terms (or omission for the extremes).  One defines $\cynerve H$ as the simplicial space whose space of $p$-simplices is given by $(p+1)$-fold product of copies of $H$, with the face map defined similarly, although the last face map is given by first cyclic permuation and then multiplying the first two terms. 
Notice there is a simplicial map $\pi:\cynerve H \to \nerve H$ given on simplices 
by projection away from the first term, which induces a composite
\[S^1 \times \rcynerve{H} \to \rcynerve{H} \to \classify {H},
\]
where the first map is given by the $S^1$-action (see \cite[page 472]{BHM}).  The adjoint map
$\rcynerve{H} \to \Lambda \classify{H}$ is then a homotopy equivalence when $H$ is group-like.

See Goodwillie \cite[I.1.8]{Goo} (or 
Bousfield and Friedlander \cite[page 311]{BF}) for the following.  Let $F$ represent the associated free monoid functor, given by realizing the simplicial resolution associated to the adjunction between topological monoids and spaces, and let $H$ represent the associated topological group functor.

\begin{lemma}
	\label{makeI}
For an FSP $L$, there is a weak natural transformation $I$
\diagramit{
\classify{\glk{L}} & \ar[l]_{\simeq} \classify{F\glk{L}} \ar[r]^{\simeq} & \classify{H\glk{L}} \ar[d]^{i} \\
	\rcynerve{\glk{L}} & \ar[l]_{\simeq} \rcynerve{F\glk{L}} \ar[r]^{\simeq} & \rcynerve{H\glk{L}} 
}
which acts as a section, in the homotopy category of spaces,
 of the natural projection $\pi:\rcynerve{\glk{L}} \to \classify{\glk{L}}$.
\end{lemma}

  Moving on toward defining the trace, there is a straightforward map 
\[ M:\rcynerve{\glk{L}} \to \thh{\mk{L}}
\] 
induced by (including and) multiplying at the level of simplices. 
Next we have the Morita equivalence map following \cite[(3.9.3) on page 480]{BHM}
as always with $r=1$ since we are not after $C_r$-equivariant results.
There is also a nice description of this map
due to Schlichtkrull \cite{Sch}, using a slightly different definition of $\THH$.

\begin{lemma}
	There is a weak natural transformation of spaces for any FSP $L$,
\diagramit{ 
\thh{\mk{L}} & \ar[l] {| | X_{\bullet, \bullet} (1) | |} \ar[r] & \thh{L}, 
}
	where both natural transformations are natural weak equivalences.
\end{lemma}

Hence, the weak composite 
\[\classify{\glk{L}} \to \rcynerve{\glk{L}} \to \thh{\mk{L}} \to  \thh{L}
\]
gives a weak natural transformation, as space valued functors from FSPs.  Of course,
we then promote to a weak natural transformation of spectra by noting all of the above
is compatible with $\Gamma$-space structures, as described in \cite[\S 4]{BHM}.
Choosing the FSP $L_{\gofx}$ to be multiplication by the (pointed version of the)
topological group $\gofx$ 
(coming from the Kan loop group functor, see subsection \ref{prodtrick}), we then
have a composite weak natural transformation as functors from spaces (after careful choices
of a basepoint for every component) to spectra
\[\kof{L_{\gofx}} \to t\thh{\mk{L_{\gofx}}} \to t\thh{L_{\gofx}} .
\]

As described above, by $\gofx$ a group-like topological group and a loop group for $X$,
there is a natural weak equivalence 
\[\rcynerve{\gofx} \to \Lambda \classify{\gofx} \to \Lambda X.
\]
See \cite[Prop. 3.7]{BHM} (with $r=1$) for details, keeping in mind that we again promote from spaces to spectra by compatibility with the $\Gamma$-space structures.

\begin{lemma}
There is a weak natural transformation of spectra
\[\thh{L_{\gofx}} \stackrel{\sim}{\to}  \qof{\Lambda \classify{\gofx}_+}
	\stackrel{\sim}{\to} \qof{\Lambda X_+}
\]
where each component map is a natural weak equivalence.  
\end{lemma}

Combining the results of this
subsection and the last then yields the following:

\begin{prop}
	\label{pretrace}
There is a weak natural transformation of functors from spaces 
(after appropriate choices of basepoints for components) to spectra
\diagramit{
\aof{X} \ar[r] &\kof{L_{\gofx}} \ar[r] & t\thh{\mk{L_{\gofx}}} \ar[r] & t\thh{L_{\gofx}} \ar[d] \\
&&  \qof{X_+} & \qof{\Lambda X_+} \ar[l]^{eval_1}   ,
}
which we refer to as the evaluation of the (extended) B\"okstedt trace and denote $\trace$.
\end{prop}

Now we can define the relative trace map as the relative weak natural transformation associated to 
$\trace$, using Lemma \ref{pairnat}.  (See Remark \ref{pairnatrem}.)

\begin{cor}
	\label{tracedef}
	There is a relative weak natural transformation, called the relative trace map,
	 $\Psi(\trace):\aof{X,Y} \to \qof{X,Y}$. 
\end{cor}

Once again, this will be natural only with respect to maps of CW pairs.

\subsection{The Candidate for the Transfer}

Since the axioms of Becker and Schultz work with a transfer of pairs, we have introduced relative versions of
the assembly, $\Aof$-theory transfer, and the evaluation of the B\"okstedt trace map.
We now have three weak natural transformations given a
perfect fibration of pairs $(p_X,p_Y):(E_X,E_Y) \to (X,Y)$ and we take their (weak) composite to form our expected weak 
map $\tof{p,E_X,E_Y}$:
\diagramit{
{\qof{X,Y}} \ar[r]^-{\Psi(\alpha)} & {\aof{X,Y}} \ar[rr]^-{\atrans{p_X,p_Y}} &&  
{\aof{E_X,E_Y}} \ar[r]^-{\Psi(\trace)} & {\qof{E_X,E_Y}}
}
which we abbreviate as $\Tof$.

For the convenience of the reader, we will next state the Becker-Schultz axioms and their main theorem in the current language.  The statement of the additivity axiom requires reference to sums of maps, which we will always perform using the $H$-space operation from the last loop space operation being performed in the usual $\colim \Omega^n \Sigma^n E$.  Our difference operation will then come from reversing the direction of the loop in the second coordinate.  These choices become irrelevant in the homotopy category, hence are only required to make our statements well-defined with respect
to weak maps.

\begin{defn} 
\label{BSaxioms} 
	A transfer is a function which assigns to each perfect fibration of pairs 
	$p:(E_X,E_Y) \to (X,Y)$ a weak map 
	$\tof{p,E_X,E_Y}:\qof{X,Y} \to \qof{E_X,E_Y}$ satisfying:
\begin{itemize}
\item \textbf{Naturality}: If $f$ is a fiberwise homotopy equivalence covering the map of
	CW pairs $g:(W,Z) \to (X,Y)$, then $\tof{p_1,E_X,E_Y} \circ \qof{g}$ and
	$\qof{f} \circ \tof{p_2,E_W,E_Z}$ are weakly equivalent (as weak maps, see \ref{weakmaps}(c)).
\item \textbf{Normalization}: $\tof{1,X,Y}$ is weakly equivalent to the identity (as a weak map).

\item \textbf{Multiplicativity}: 	
Given two perfect fibrations of pairs $p_1:(E_X,E_Y) \to (X,Y)$ and $p_2:(E_W,E_Z) \to (W,Z)$,
	the diagram 
\diagramit{
{\qof{X,Y} \Smash \qof{W,Z}} \ar[d]_{\tof{p_1} \Smash \tof{p_2}} \ar[r]^-{\mu} 
    & {\qof{X \times W,X \times Z \cup Y \times W}} 
    \ar[d]^{\tof{p_1 \times p_2}} \\
{\qof{E_X,E_Y} \Smash \qof{E_W,E_Z}} \ar[r]^-{\mu} 
    & {\qof{E_X \times E_W,E_X \times E_Z \cup E_Y \times E_W}} .
}
	commutes up to a natural homotopy in the category of spectra.
		
\item \textbf{Strong Additivity}:  Assume
\diagramit{
{E^1} \ar[dr]_{p_1} & {E^0} \ar[l]_{i} \ar[r]^{j} \ar[d]^{p_0} 
    & {E^2} \ar[dl]^{p_2} \\
& {X}
}
is a diagram of perfect fibrations, with
$p:E_X \to X$ the pushout (perfect) fibration and $i$ a cofibration,
while $(X,Y)$ is a CW-pair.  
For $n=0,1,2$, let $k^n:(E^n_X,E^n_Y) \to (E_X,E_Y)$ 
denote the relevant inclusion 
and $p_n:(E^n_X,E^n_{Y}) \to (X,Y)$ the
perfect fibration of pairs.  Then $\tof{p}$ and
$ \qof{{k}^1} \tof{p_1}  + \qof{{k}^2}\tof{p_2} -  \qof{{k}^0}\tof{p_0}$
are weakly equivalent (as weak maps).

\end{itemize}
\end{defn}

\begin{thm} 
	\label{BSthm}
	Any transfer (in the sense of Def. \ref{BSaxioms}) is weakly equivalent
	(as a weak map, see \ref{weakmaps}(c)) to the Becker-Gottlieb transfer when evaluated at
	any compact ANR fibration over a finite base.
\end{thm}

\begin{proof}
	This is merely a restatement of the Main Theorem of \cite{BecSch}
	into the language of weak maps.  Since a transfer as in \ref{BSaxioms} induces a bundle
	transfer in the notation of \cite{BecSch} satisfying axioms I-III and IV$^+$, their result 
	implies the induced map in the homotopy category is the Becker-Gottlieb transfer.
	As the construction of the Becker-Gottlieb transfer in either \cite{BecGot}
	or \cite{BecGot2} does produce a weak map, it then
	suffices to recall our convention that weak maps be called weakly
	equivalent when they induce the same map in the homotopy category.
\end{proof}

We can now state our main result:

\begin{thm}
	\label{main}
	Suppose $X$ is a finite CW complex and $p:E \to X$ is a compact ANR fibration.
	Then the evaluation of the B\"okstedt trace of the algebraic $\Kof$-theory transfer of $p$ is	weakly equivalent (as a weak map, see \ref{weakmaps}(c)) to the
Becker-Gottlieb transfer of $p$, that is $\trace \circ \ktrans{p} = \bgtrans{p}$ in the 
homotopy category of spectra.
\end{thm}

\begin{proof}
As indicated above, our method of proof is to verify the axioms of Becker and Schultz
for $\tof{p}$ (assuming only that $p$ is a perfect fibration with finite base).
We verify the naturality and normalization axioms 
in Props \ref{naturality} and \ref{normalization},
the multiplicativity axiom in Prop \ref{prodtransfer}
and the strong additivity axiom in Prop \ref{additivity}.
\end{proof}

\section{Multiplicativity}

The point here is to show that all three (relative) components of the candidate transfer
$\Tof$ are compatible with the external pairing on relative $\Aof$-theory.  

\subsection{Multiplicativity of the $\Aof$-theory transfer}

We begin by introducing relative pairings and our desired notion of multiplicative, 
keeping in mind the material from 
subsection \ref{weaknat} on induced relative functors, denoted there with a $\Psi$.

\begin{defn}
	\label{multweak}
	When each intermediate spectrum-valued functor in a 
	weak natural transformation (in the sense of subsection 2.6) is equipped
	with an external pairing, such that each diagram of the form
\diagramit{
{F(X) \Smash F(W)} \ar[d]_{\zeta \Smash \zeta} \ar[r]^{\mu_F} & {F(X \times W)} \ar[d]^{\zeta} \\ 
{G(X) \Smash G(W)} \ar[r]^{\mu_G} & {G(X \times W)}
}
	commutes up to a homotopy natural in each variable, it
	will be referred to as a multiplicative weak natural transformation.
\end{defn}

\begin{lemma}
	\label{relmult}
Let $F$ be a spectrum-valued functor of spaces, equipped with an external pairing that is
natural up to a natural homotopy.  Then the external product for the functor $F$ induces a 
relative external product
\[
\mu_{rel}:F(X,Y) \wedge F(W,Z) \rightarrow F(X \times W, X\times Z \cup _{Y \times Z}Y \times W)
\]
which is also natural up to a natural homotopy.
\end{lemma} 

\begin{proof}
	First, notice the left hand side is naturally homotopy equivalent
	to the homotopy pushout of the diagram
\def\objectstyle{\scriptstyle}
\diagramit{
\hocofib(i:F(Y) \Smash F(Z) \to F(X) \Smash F(W)) \ar[r] \ar[d] 
	& \hocofib(j:F(Y) \Smash F(W) \to F(X) \Smash F(W)) \ar[d] \\
\hocofib(k:F(X) \Smash F(Z) \to F(X) \Smash F(W)) \ar[r] & P . 
}
	Now notice 
\[F(X \times W, X\times Z \cup _{Y \times Z}Y \times W)
\]	
	is also defined as a homotopy cofiber. Furthermore, from all three homotopy cofibers 
	in the homotopy pushout diagram above,
	the maps from the targets of $i$,$j$, and $k$
\[\mu: F(X) \Smash F(W) \to F(X \times W) 
\]
 	are the same.
	Thus, it will suffice to verify that the three composite maps $\mu \circ j$,
	$\mu \circ i$ and $\mu \circ k$ each factor up to homotopy
	through the map
\[F(X\times Z \cup _{Y \times Z}Y \times W) \to F(X \times W).
\]	
	This follows from naturality of the homotopy for the external product of $F$
	 and the commutative pushout diagram which defines the union.
\end{proof}

The key application is the relative pairing in $\Aof$-theory that follows.

\begin{cor}
  \label{relexternal}
  The external smash product of retractive spaces induces an 
  external pairing on relative $\Aof$-theory
\[\mu:\aof{X,Y} \Smash \aof{W,Z} \to \apof{X}{Y}{W}{Z}.
\] 
\end{cor}

\begin{proof}
	By Lemma \ref{relmult}, it will suffice to show the diagram of Waldhausen categories 
	and (bi-)exact functors
\diagramit{
\retof[Y] \times \retof[Z] \ar[r]^-{\Smash_{Y \times Z}} \ar[d]_{i_* \times i'_*} &
		\retof[Y \times Z] \ar[d]^{(i \times i')_*}  \\
\retof \times \retof[W]  \ar[r]^-{\Smash_{X \times W}} & \retof[X \times W]		
}
	commutes up to a unique isomorphism.
	
	Recall product with a fixed space preserves pushouts (as we work with compactly
	generated spaces), and the universal property of pushouts implies there are unique 
	isomorphisms
	$i_*Y \to X$, $i'_*Z \to W$, and $(i \times i')_*(Y \times Z) \to X \times W$.
	Hence, for any pair $U \in \retof[Y]$ and $V \in \retof[Z]$, there is also a unique 
	isomorphism
\[i_*U\hspace{.03in} _X\hspace{-.05in}\vee_{W} i'_*V \to (i \times i')_* (U\hspace{.03in} _Y\hspace{-.05in}\vee_{Z} V)
\]
	by inspection of the defining pushouts (see subsection \ref{retrdef}).  Now by inspection of
	the defining pushout diagrams, one has unique isomorphisms in all three input slots, 
	hence the required unique isomorphism
\[i_*U\hspace{.03in} _X\hspace{-.05in}\Smash_{W} i'_*V \to (i \times i')_* (U\hspace{.03in} _Y\hspace{-.05in}\Smash_{Z} V).
\]
\end{proof}

\begin{lemma}
	\label{transmult}
	Given two perfect fibrations of pairs $p_1:(E_X,E_Y) \to (X,Y)$ and $p_2:(E_W,E_Z) \to (W,Z)$,
	the diagram 
\diagramit{
{\aof{X,Y} \Smash \aof{W,Z}} \ar[d]_{\atrans{p_1} \Smash \atrans{p_2}} \ar[r]^-{\mu} 
    & {\aof{X \times W,X \times Z \cup Y \times W}} 
    \ar[d]^{\atrans{p_1 \times p_2}} \\
{\aof{E_X,E_Y} \Smash \aof{E_W,E_Z}} \ar[r]^-{\mu} 
    & {\aof{E_X \times E_W,E_X \times E_Z \cup E_Y \times E_W}} .
}
	commutes up to a natural homotopy in the category of spectra.
\end{lemma}

\begin{proof}
  Since $\Omega|Nw\sdot(?)|$ 
  preserves products up to isomorphism, it would suffice to verify that 
  the following diagram of
  (bi-)exact functors commutes up to a unique natural isomorphism 
  (which as a consequence will be compatible with
  restriction to subspaces):
\diagramit{
{\retof \times \retof[W]} \ar[d]_{p_1^* \times p_2^*} \ar[rr]^-{_X \hspace{-.02in}\Smash_{W}} && 
  {\retof[X \times W]} \ar[d]^{(p_1\times p_2)^*} \\
{\retof[E_X] \times \retof[E_W]} \ar[rr]_-{_{E_X}\hspace{-.02in}\Smash_{E_W}} && 
  {\retof[E_X \times E_W]}.
}
  Choose $U \in \retof$ and $V \in \retof[W]$.	
  Since $X \rightarrowtail U$ and $W \rightarrowtail V$ are
  cofibrations, we see 
\[\pof{U}{X}{V}{W} \rightarrowtail U \times V
\]
  is also a cofibration.  Hence, we again use Lemma \ref{wolfgang} to
  see there is a unique natural isomorphism in $\retof[E_X \times E_W]$
\[p_1^*(U)\hspace{.03in} _{E_X}\hspace{-.05in} \Smash_{E_W} p_2^*(V) \cong 
    (p_1 \times p_2)^*(U\hspace{.03in} _X \hspace{-.05in} \Smash_{W} V).
\]
\end{proof}

\subsection{Multiplicative Relative Weak Natural Transformations}

\begin{defn}
	\label{multweakrel}
	An induced relative weak natural transformation $\Psi(\zeta)$ will be called
	multiplicative if for each component natural transformation
	$\zeta:F \to G$ the source and target are equipped with external pairings
	and the following diagram commutes 
	up to a homotopy which is natural in each variable:
\diagramit{
F(X,Y) \wedge F(W,Z) \ar[r]^-{\mu_F} \ar[d]_{\Psi(\zeta) \Smash \Psi(\zeta)} & 
	F(X \times W, X\times Z \cup _{Y \times Z}Y \times W) \ar[d]^{\Psi(\zeta)} \\
G(X,Y) \wedge G(W,Z) \ar[r]^-{\mu_G} &
G(X \times W, X\times Z \cup _{Y \times Z}Y \times W)
}
\end{defn}

The key general result in this direction is the following:

\begin{lemma}
	\label{yuck}
	Suppose $\zeta:F \to G$ is a multiplicative weak natural transformation (in the
	sense of Def. \ref{multweak}).  Then $\Psi(\zeta)$,
	the relative weak natural transformation it induces, is multiplicative (in the
	sense of Def. \ref{multweakrel}).
\end{lemma}

\begin{proof}
	As in the proof of Lemma \ref{relmult}, recall $F(X,Y) \Smash F(W,Z)$
	is naturally homotopy equivalent
	to a homotopy pushout, as is $G(X,Y) \Smash G(W,Z)$.
	Thus, it will suffice to show the two maps 
\[ \hocofib(F(Y) \Smash F(W) \to F(X) \Smash F(W)) \to G(X \times W,\pof{X}{Y}{W}{Z}) 
\]
	are homotopic in a natural way, and similarly for $F(X) \Smash F(Z)$.  
	Notice that each of these maps factors through
	$G(X \times W,Y \times W)$,
	so it will suffice to see they are naturally
	homotopic up to that point.  Since naturality of this
	homotopy will imply compatibility with the ``restriction" to 
\[\hocofib(F(Y) \Smash F(Z) \to F(X) \Smash F(W)),
\]
	this will suffice to complete the proof.

	Now recall that a map from one homotopy cofiber of this form to another is completely determined
	by the corresponding null-homotopy of the composite 
\[F(Y) \Smash F(W) \to G(Y \times W) \to G(X \times W).
\]
	Thus, it will suffice to see the two null-homotopies are themselves naturally homotopic.
	In the case of the first portion of $\mu_G \circ (\zeta \Smash \zeta)$, 
	the null homotopy will arise by applying the reduced cone operator 
	(smashing the 
	composite map	 with the unit interval, using one end as basepoint)	to 
\[F(Y) \Smash F(W) \to G(Y) \Smash G(W) \to G(Y \times W)
\]
	and 
	then following with the chosen null homotopy of $G(Y \times W) \to G(X \times W)$.  For the 
	first portion of $\zeta \circ \mu_F$, the construction is similar using the composite
\[F(Y) \Smash F(W) \to F(Y \times W) \to G(Y \times W).
\]
	Since these two composites are assumed to be naturally homotopic by the assumption that
	$\zeta$ is a multiplicative weak natural transformation, we may extend this natural
	homotopy over the cones to show the two induced null homotopies of 
\[F(Y) \Smash F(W) \to G(Y \times W) \to G(X \times W)
\]
	are naturally homotopic, which suffices to complete the proof as indicated above.
\end{proof}

\begin{prop}
	\label{multdeal}
	Both the assembly map and the evaluation of the B\"okstedt trace map are multiplicative
	weak natural transformations of spectra.
\end{prop}

\begin{proof}
	This is a tedious business of inspecting each commutative square involved in the 
	construction of the assembly map (see subsection \ref{makerelass}) and in the construction
	of the transfer map (see subsection \ref{makereltrace}).  
	
	For the assembly map, the first three functors involved all have external pairing maps
	induced from natural homotopy equivalences on the variables together with a ring
	spectrum product, which makes the condition simple to verify.
	The remaining two portions come from inspecting the external pairing on $\aperof{?}$,
	which comes from the external smash product of retractive spaces over simplices.  
	Taken together, this implies the assembly map is a multiplicative weak natural transformation.
		
	For the trace map,  we begin by verifying that the change of $\Aof$-theory models from the 		
	retractive spaces model to the FSP model is multiplicative:  
	
	Recall the transition (of Waldhausen's \cite{Wal3}) from the retractive spaces
	over $X$ model to the free (pointed) $\gofx$-spaces model takes a retractive space 
	$W$ over $X$ to $p_1^*(W)/P_X$, where $P_X$ is a universal $\gofx$ bundle over $X$. 
	In the retractive spaces model the external smash product induces the external pairing, 
	while the ordinary smash product induces the external pairing in the free $\gofx$-spaces model.
	Since $P_{X \times Y}$ 
	is homeomorphic to $P_X \times P_Y$, taking the external smash product of $p_1^*(W)$ 
	and $p_2^*(V)$ and then collapsing
	$P_{X \times Y}$ will be homeomorphic to taking the smash product of collapsing $P_X$ in 
	$p_1^*(W)$ and $P_Y$ in $p_1^*(V)$ separately.  This will give a natural homeomorphism
	at the level of models, hence a natural homotopy at the level of $\Aof$-theory spectra,
	since these transitions are all given by exact functors.	
	
	For the remaining steps in the passage from the free $\gofx$-spaces model to Waldhausen's
	version of the FSP model,
	at each step the external pairing is induced by the ordinary smash product of pointed
	spaces.  Thus, an inspection of the chain of weak homotopy equivalences in a 
	proof of Waldhausen 
	\cite[Theorem 2.2.1, pages 386 through 388]{Wal3} 
	leads to the conclusion that the transitions between
	various models are all compatible (up to a homotopy natural in each variable)
	with external pairings in the homotopy category of spectra.
	
	For the passage from Waldhausen's version of the FSP model to $\kof{L_{\gofx}}$,
	consider the following.  The naturality of $\gofx$
	yields a natural map $\gofx[X \times Y] \to \gofx \times \gofx[Y]$ which is an
	isomorphism by inspection \cite{WalK} and induces the external pairings in question
	(see \cite[(1.3)]{Wal2}).   For FSPs, this gives a pairing (in the sense of \cite{BHM}) 
\[L_{\gofx} \Smash L_{\gofx[Y]} \to L_{\gofx \times \gofx[Y]}
\]
	whose target is canonically identified to $L_{\gofx[X \times Y]}$ by the isomorphism
	of topological groups.
	Careful inspection of the construction at the heart of the proof of \ref{FSPmodel} then verifies the
	resulting external pairing is natural on the level of spaces.  In order to promote to spectra,
	as we will continue to do below, one then verifies all constructions are compatible with the
	$\Gamma$-space structures defined in \cite[\S 2.3.4]{DGM} (or later, in \cite[\S 4]{BHM}).
	
	Next B\"okstedt, Hsiang, and Madsen \cite[page 503]{BHM} point out that 
	the B\"okstedt trace map from the FSP model for $\Aof$-theory
	to $t\THH$ is multiplicative in the homotopy category of spectra.  Since their pairings are 
	all induced by the pairing of FSPs indicated above,	it is simply a tedious verification to 
	see the relevant homotopy can be chosen to be natural in each variable.
		
	For the passage from
	$\THH$ to the stable free loop space, notice the interlacing of spheres built into the external 
	product on $\THH$ reduces this portion to verifying that the map
	$|\cynerve{\Gamma_X}| \to \Lambda |\nerve{\Gamma_X}|$ is multiplicative as in 
	\cite[page 472]{BHM}.  To see the remaining portion,
	notice the standard $S^1$ action on
	a product of cyclic nerves is the diagonal action, where $S^1$ acts separately on each 
	factor, and this is compatible with the usual action on the cyclic nerve of a product.
	Thus, taking the adjoint yields a product of maps as above, compatible with the 
	adjoint map for the product FSP.  Thus, the map from $\THH$ to the stable free loop space
	is multiplicative up to a natural homotopy.
	
	 Finally, evaluation at the unit of $\mathbb C$ will be multiplicative on the level of spaces
	 by construction, and the stabilization preserves smash products up to natural isomorphism.
	 Thus, the long composite is, in fact, a multiplicative weak natural transformation.
\end{proof}

Combining the last two results now gives us what we will need below.

\begin{prop}
 	\label{prodtransfer}
 	Given two perfect fibrations of pairs $p_1:(E_X,E_Y) \to (X,Y)$ and $p_2:(E_W,E_Z) \to (W,Z)$,
	the diagram 
\diagramit{
{\qof{X,Y} \Smash \qof{W,Z}} \ar[d]_{\tof{p_1} \Smash \tof{p_2}} \ar[r]^-{\mu} 
    & {\qof{X \times W,X \times Z \cup Y \times W}} 
    \ar[d]^{\tof{p_1 \times p_2}} \\
{\qof{E_X,E_Y} \Smash \qof{E_W,E_Z}} \ar[r]^-{\mu} 
    & {\qof{E_X \times E_W,E_X \times E_Z \cup E_Y \times E_W}} .
}
	commutes up to a natural homotopy in the category of spectra.
\end{prop}

\begin{proof}
	Notice $\Psi(\alpha)$ and $\Psi(\trace)$ are multiplicative relative weak natural
	transformations by Lemma \ref{yuck} and Prop \ref{multdeal}, while the multiplicativity of
	the relative $\Aof$-theory transfer is verified in Lemma \ref{transmult}.
\end{proof}

\section{Naturality and Normalization}

This section contains two propositions, dealing with each of the axioms listed in the section title for our candidate transfer $\Tof$.  With the technical preliminaries handled in section 2, each of these is straightforward, unlike the remaining strong additivity axiom which is considered in the last section. 

First, we establish the naturality property of $\Tof$.  

\begin{prop}
  \label{naturality}
Suppose $f:E_W \to E_X$ is a fiberwise homotopy equivalence making
the diagram
\diagramit{
{(E_W,E_Z)} \ar[r]^{f} \ar[d]_{p_2} & {(E_X,E_Y)} \ar[d]^{p_1} \\
{(W,Z)} \ar[r]_{g} & {(X,Y)}
}
commute.
Then the diagram
\diagramit{
{\qof{W,Z}} \ar[r]^{\qof{g}} \ar[d]_{\tof{p_2,E_W,E_Z}}  
         & {\qof{X,Y}} \ar[d]^{\tof{p_1,E_X,E_Y}} \\
{\qof{E_W,E_Z}} \ar[r]_{\qof{f}} & {\qof{E_X,E_Y}}
}
commutes in the homotopy category of
spectra (or equivalently, the composites are weakly equivalent as weak maps).
\end{prop}

\begin{proof}
This breaks up into a series of 
commutative diagrams in the homotopy category of spectra.
Commutativity of 
\diagramit{
{\qof{W,Z}} \ar[r]^{\qof{g}} \ar[d]_{\Psi(\alpha)}  
         & {\qof{X,Y}} \ar[d]^{\Psi(\alpha)} \\
{\aof{W,Z}} \ar[r]_{\aof{g}} & {\aof{X,Y}}
}
up to homotopy is the naturality of the relative assembly map from Cor \ref{relassem},
since $g$ is a map of pairs.
Then homotopy commutativity of the diagram  
\diagramit{
 {\aof{W,Z}} \ar[r]^{\aof{g}} 
\ar[d]_{\atrans{p_2}}
         & {\aof{X,Y}} \ar[d]^{\atrans{p_1}} \\
{\aof{E_W,E_Z}} 
\ar[r]_{\aof{f}}& {\aof{E_X,E_Y}} 
}
follows by Lemma \ref{transdef}.
Finally, we have the homotopy commutative diagram
\diagramit{
{\aof{E_W,E_Z}}  \ar[r]^{\aof{f}} \ar[d]_{\Psi(\trace)} 
    & {\aof{E_X,E_Y}} \ar[d]^{\Psi(\trace)} \\
{\qof{E_W,E_Z}}  \ar[r]^{\qof{f}}& {\qof{E_X,E_Y}}
}
given by the naturality of the relative trace map from Cor \ref{tracedef}, 
since $f$ is also a map of pairs.
\end{proof}

Now we have the normalization axiom for $\Tof$.

\begin{prop}
  \label{normalization}
  For the identity map $1:(X,Y) \to (X,Y)$, the
  weak map
\[\tof{1,X,Y}:\qof{X,Y} \to \qof{X,Y}.
\]
  is weakly equivalent to the identity.
\end{prop}

\begin{proof}
Since the relative $\Aof$-theory transfer of the identity fibration is clearly 
equivalent to the identity by construction, it suffices to see that 
\[\qof{X,Y} \overset{\alpha}{\to} \aof{X,Y} 
    \overset{\trace}{\to} \qof{X,Y}
\]
is the identity up to a natural weak equivalence.

However, in the absolute case Waldhausen shows  \cite[Theorem 5.1]{Wal2}
\[\qof{X_+} \overset{\alpha}{\to} \aof{X} 
    \overset{\trace}{\to} \qof{X_+}
\]
is naturally weakly equivalent to the
identity on $\qof{X_+}$
(as a map of spectra by his Remark 5.3), which implies the same is
true with $Y$ in place of $X$ everywhere.  Hence 
naturality of this weak equivalence implies 
$\Psi(\trace \alpha) \sim \Psi(\trace)\Psi(\alpha) \sim \tof{1,X,Y}$ 
is naturally weakly equivalent to the identity as well.
\end{proof}

\section{Strong Additivity}

We now move on to the most subtle of the axioms for our composite, the strong additivity 
axiom, verified in Proposition \ref{additivity}.  One would like to simply appeal to the additivity of the evaluation of 
the B\"okstedt trace and additivity of the algebraic $\Kof$-theory transfer map.
Unfortunately, the homotopy constructed by the first author \cite{Do} to verify the additivity of the algebraic $\Kof$-theory transfer need not be natural, because of its reliance upon Waldhausen's Additivity Theorem \cite[Thm 1.4.2]{Wal3} which makes no claim of naturality.  Thus, to prove 
additivity for the algebraic $\Kof$-theory transfer map in the relative case, we must verify the 
existence of compatible homotopies for the sub-fibration and the original fibration.  Throughout 
this section we will be working with retractive spaces, with notation established in section 2.

As referred to in subsection \ref{bivardef}, Williams \cite{Wil} uses the homotopy (parametrized) Euler characteristic to give lifts of transfer maps for the algebraic $\Kof$-theory of spaces to his bivariant 
$\Aof$-theory.  Thus, our required homotopies may be produced by constructing paths in this bivariant $\Aof$-theory space between points associated to homotopy Euler characteristics, which is the technique used by the first author in \cite{Do}.  We will expand upon this technique, using a different map out of the bivariant $\Aof$-theory, to produce our homotopy for the sub-fibration so that it will be compatible with that for the original fibration.

We begin by introducing some notation we will use throughout this section.

We will assume
\diagramit{
{E^1} \ar[dr]_{p_1} & {E^0} \ar[l]_{i} \ar[r]^{j} \ar[d]^{p_0} 
    & {E^2} \ar[dl]^{p_2} \\
& {X}
}
is a diagram of perfect fibrations, with
$p:E_X \to X$ the pushout (perfect) fibration and $i$ a cofibration,
while $(X,Y)$ is a CW-pair.  
For $n=0,1,2$, let $k^n:(E^n_X,E^n_Y) \to (E_X,E_Y)$ 
denote the relevant inclusion 
and $p_n:(E^n_X,E^n_{Y}) \to (X,Y)$ the
perfect fibration of pairs.

The heart of our additivity result is the  $\Aof$-theory
Transfer Additivity Theorem of the first author \cite{Do}, which suggests the following definition.  
Recall the sum operation used here comes from the coproduct operation 
\[(W,Z) \mapsto W \sqcup Z
\]
in the category of retractive spaces, while the difference operation is instead induced by 
\[(W,Z) \mapsto W \sqcup \Sigma Z 
\]
where $\Sigma$ here indicates the fiberwise suspension in the category of retractive spaces.

\begin{defn}
	\label{decomp}
	Given a decomposed perfect fibration $p_X$ as above, we will let
	$p_{dec(X)}^*:\aof{X} \to \aof{E_X}$ denote the sum (indicated above) of maps 
\[\aof{k^1} p_1^* + \aof{k^2} p_2^* - \aof{k^0} p_0^*.
\]
\end{defn}

\begin{thm}[Transfer Additivity Theorem of Dorabiala \cite{Do}]
  \label{wojtek}
  For a pushout fibration $p_X$ as above, 
  there is a homotopy in the category of spectra
\[p_X^* \simeq \aof{k^1} p_1^* + \aof{k^2} p_2^* 
   - \aof{k^0} p_0^*=p_{dec(X)}^*.
\]
\end{thm}

Given a CW-pair $(X,Y)$ as above, we would like to define the decomposed transfer of the perfect fibration of pairs.  Thus, we first note that there is a natural isomorphism 
(at the level of retractive spaces) between the composites in the diagram
\diagramit{
{\retof[Y]} \ar[r] \ar[d]_{p_j^*} & {\retof} \ar[d]^{p_j^*} \\
{\retof[E^j_Y]} \ar[r] & {\retof[E^j]} .
}
As a consequence, we may define the decomposed transfer of the perfect fibration of pairs
as the homotopy class of maps induced (as in Def. \ref{makepairs}) 
from the corresponding homotopy in the diagram of spectra below.
\diagramit{
{\aof Y} \ar[r] \ar[d]_{p_{dec(Y)}^*}  & 
	{\aof {X}} \ar[r] \ar[d]_{p_{dec(X)}^*} &
	{\aof {X,Y}} \ar@{-->}[d]_{\atransdec{p_{X,Y}}}  \\
{\aof{E_Y}} \ar[r] & {\aof{E_X}} \ar[r] & {\aof{E_X,E_Y}}
}

We will verify that the homotopy Euler characteristic, considered as a point in the zero space of the
bivariant $\Aof$-theory spectrum, is a lift of $p_X^* \alpha$ (see Lemma \ref{samecomp}) over a certain map of spectra
\[\aofp[X] \to \mapsp{X,\aof{E_X}}.
\]
Similarly, one can define (see \cite[page 260]{Do}) a point $\chi(p_{dec(X)})$ in $\Omega^\infty \aofp[X]$ which is a lift
of $p_{dec(X)}^* \alpha$ over the displayed map.
These are the same points in $\aofp[X]$ shown to be connected by a (not neccesarily natural) path in 
\cite{Do}, and we will label that path $\gamma_X$.
Williams's bivariant theory then comes equipped with a restriction map $Res:\aofp[X] \to \aofp[Y]$ coming from the pullback diagram defining a sub-fibration.  However, it is not clear that 
$Res(\gamma_X)=\gamma_Y$, although there is a natural path 
$\beta_{X,Y}$ from $Res(\chi(p_{X}))$ to $\chi(p_{Y})$
(since they correspond to retractive spaces which are naturally homeomorphic).  
One also has a similar path 
$\beta_{dec(X,Y)}$ from $Res(\chi(p_{dec(X)}))$ to $\chi(p_{dec(Y)})$, again corresponding to a
natural homeomorphism of retractive spaces.  Thus, rather than choosing 
$\gamma_Y$, we may choose 
\[\omega_Y=\beta_{X,Y}^{op} * Res(\gamma_X) * \beta_{dec(X,Y)}
\]
which satisfies $Res(\gamma_X)$ is homotopic to $\omega_Y$.  At the end of this section, we will produce a homotopy commutative diagram of spectra (see Lemma \ref{techone})
\diagramit{
{\aofp[X]} \ar[d]_{\rho^{E_X}\nu^{E_X}} \ar[rr]^{Res} 
   && {\aofp[Y]} \ar[d]^{\rho^{E_Y}\nu^{E_Y}} \\
{\mapsp{X,\aof{E_X}} } \ar[dr]_{\epsilon}  &&  
	   {\mapsp{Y,\aof{E_Y}}} \ar[dl]^{\delta} \\
& {\mapsp{Y,\aof{E_X}}}.
}	
As a consequence, one may choose $H_X$ to be the image of $\gamma_X$ and $H_Y$ to be the image of $\omega_Y$ so $\epsilon H_X \simeq \delta H_Y$ implies the homotopies are appropriately compatible. Thereby one produces the required (induced) homotopy between the maps on homotopy cofibers, from $\atrans{p_{X,Y}}\Psi(\alpha)$ to $\atransdec{p_{X,Y}}\Psi(\alpha)$. 

We now proceed with the verification of the strong additivity property for our candidate transfer $\Tof$, based upon the technical results which will fill the remainder of this section.

\begin{prop}
  \label{additivity}
For a pushout perfect fibration of pairs $p$ as above, there is a weak
equivalence (of weak maps)
\[
\tof{p} \sim \qof{{k}^1} \tof{p_1}  + \qof{{k}^2}\tof{p_2} -  \qof{{k}^0}\tof{p_0},
\]
that is, the composites agree in the homotopy category of spectra.
\end{prop}

\begin{proof}
	This is a consequence of Proposition \ref{naturaladd}
	and natural additivity of the assembly map (up to 
	a natural homotopy in the category of spectra),
	which combine to say there is a weak equivalence (of weak maps)
\[
\atrans{p}\Psi(\alpha)  \sim 
    \aof{{k}^1} \atrans{p_1}\Psi(\alpha) + 
   \aof{{k}^2}  \atrans{p_2}\Psi(\alpha) - 
   \aof{{k}^0} \atrans{p_0}\Psi(\alpha).
\]	 
	Then we use the fact that the trace map is a 
	natural additive map (again, up to a natural homotopy in the
	category of spectra) to see 
\begin{align*}
&\Psi(\trace)  \circ 
[\aof{{k}^1} \atrans{p_1}\Psi(\alpha) + 
   \aof{{k}^2}  \atrans{p_2}\Psi(\alpha) - 
   \aof{{k}^0} \atrans{p_0}\Psi(\alpha) ]
   \sim \\
   & \qof{{k}^1} \Psi(\trace) \atrans{p_1}\Psi(\alpha) + 
   \qof{{k}^2} \Psi(\trace) \atrans{p_2}\Psi(\alpha) - 
   \qof{{k}^0} \Psi(\trace) \atrans{p_0}\Psi(\alpha)  .
\end{align*}
\end{proof}

We next begin verifying the technical background by working toward the verification that the homotopy Euler characteristic is a lift of $p_X^* \alpha$.  Recall the bivariant $\Aof$-theory space of a fibration, and the parametrized Euler characteristic of Williams \cite{Wil}, discussed in subsection \ref{bivardef}.

\begin{defn}
	\label{bar}
	Let $\simp X$ denote the category (poset) of simplices of $X$ under inclusion.
	Then $\abof[X]{E_X} $ will denote the spectrum 
	$\hocolimdown{\sigma \in \simp X}  \aof{E_X^\sigma}$,
	where $E_X^\sigma$ is the restriction of $E_X$ over $\sigma \subset X$. 
\end{defn}

We will be carefully considering the composite map of spectra
\diagramit{
{\aofp[X]} \ar[r]^-{\nu^{E_X}} & {\holimdown{\sigma \in \simp X} \aof{E_X^\sigma}} \ar[r]^{\rho^{E_X}} & 
{\mapsp{ X, \abof[X]{E_X} }} \ar[d]^{(u^{E_X})_*} \\
&&{\mapsp{X, \aof{E_X} }} .
}
Here $\nu^{E_X}$ is a Thomason homotopy inverse limit map 
(as in \cite[page 12]{Wil}).  The map
$\rho^{E_X}$ is defined, 
using the Bousfield-Kan \cite{BK} models for homotopy (co)limits,
essentially by taking a natural transformation (on $\simp X$) from 
$\qof{|\simp X/?|_+}$ to $\aof{E_X^?}$ and sending it to the induced map on hocolims and then exploiting the
usual isomorphism of spectra
\[\specmap{\qof{|\simp X|_+},\abof[X]{E_X}} \approx \mapsp{ X, \abof[X]{E_X} }.  
\]
Finally, 
$u^{E_X}:\abof[X]{E_X} \to \aof{E_X}$ is the natural map constructed from the definitions.

With this definition, we also have a description of 
$\nu^{E_X}(\chi(p)) \in \Omega^\infty\holimdown{\sigma \in \simp X} \aof{E_X^\sigma}$
as the class which in each $\Omega^\infty\aof{E_X^\sigma}$ corresponds to the retractive space 
\[E_X^\sigma \sqcup E_X^\sigma \approx E_X^\sigma \times S^0 \leftrightarrows E_X^\sigma.
\]
A decomposed version can also be described fairly explicitly (see \cite[page 256]{Do}).
Furthermore, both of these definitions are extended to the relative case using 
Def. \ref{makepairs}.

\begin{lemma}
	\label{samecomp}
	The image of $\chi(p_X)$ under $\Omega^\infty (u^{E_X})_* \rho^{E_X} \nu^{E_X}$ is (the 
	adjoint of) the composite of the assembly and the transfer
\diagramit{
{\qof{X_+}}  \ar[r]^{\alpha} & {\aof{X}} \ar[r]^-{p_X^*} & {\aof{E_X}}.
}
	Similarly, the image of $\chi(p_{dec(X)})$ under $\Omega^\infty (u^{E_X})_* \rho^{E_X} \nu^{E_X}$	is the adjoint of the expected loop space sum
\[\aof{{k}^1}  p_1^*\alpha + 
 \aof{{k}^2} p_2^*\alpha - 
   \aof{{k}^0}  p_0^*\alpha
\]
	 of the composites of the 
	assembly and transfer maps for the pieces of the decomposition 
	(pushed into the total space). 
\end{lemma}

\begin{proof}
We begin with establishing the following homotopy commutative diagram
\diagramit{
{\qof{ X_+}} \ar[r] \ar[dr] & {\abof[X]{X}} \ar[d] \ar[r] & {\aof{X}} \ar[d]^{p_X^*} \\
& {\abof[X]{E_X}} \ar[r] & {\aof{E_X}}
}
where the vertical maps are given by pullback over $p_X$ at the level of retractive spaces.
The square commutes up to homotopy by naturality (up to isomorphism) of the pullback 
construction of retractive spaces.  To see the triangle commutes up to homotopy, consider the 
following homotopy commutative diagram, where the vertical maps are induced by
transfers over $p_X$
\diagramit{
{\holimdown{\sigma \in \simp X}\aof{\sigma}} \ar[d]^{(p_X^*)_*} \ar[r]^{\rho_X} &  
{\mapsp{ X,\abof[X]{X}}} \ar[d]^{(p_X^*)_*} \\
{\holimdown{\sigma \in \simp X}\aof{E_X^\sigma}} \ar[r]^{\rho_{E_X}}\ar[r] 
& {\mapsp{ X,\abof[X]{E_X}}}. 
}
Here, the image of the $\nu^X(\chi(id))$ under $(p_X^*)_* \rho_X$
is the upper composite in the 
triangle diagram.
Similarly, the image of $\nu^X(\chi(id))$ around the lower path through this diagram is 
(up to homotopy) the lower map in the triangle by definition, since the transfer of $\nu^X(\chi(id))$
agrees with $\nu^{E_X}(\chi(p))$ by the explicit description of these classes given above.

Now the lower composite in the pentagon diagram represents the image of $\chi(p_X)$ under 
$\Omega^\infty (u^{E_X})_* \rho^{E_X} \nu^{E_X}$ by construction.  However, since the identity fibration is a smooth fiber bundle, Theorems 5.4 and 8.5 of Dwyer, Weiss, and Williams \cite{DWW} along with naturality of (fiberwise) assembly imply the horizontal composite across the top of the pentagon diagram is homotopic to the assembly map. Since the right vertical in the pentagon diagram is the transfer of $p_X$, this completes the proof for that claim.

The claim for $\chi(p_{dec(X)})$ follows from above by looking at each $p_i$ separately, since
the maps $\nu^{E_X}$,  $\rho^{E_X}$, and $u^{E_X}$ are all natural in the variable $E_X$.
\end{proof}

Let $H_X$ indicate the homotopy which is the image of $\gamma_X$ in $\mapsp{ X,\aof{E_X}}$, and $H_Y$ similarly the image of $\omega_Y$.  We will also use $l:E_Y \to E_X$ to indicate the inclusion on total spaces of the sub-fibration.

\begin{prop}
	\label{naturaladd}
	Given a CW-pair $(X,Y)$ as above, 
	there is a natural (in pairs) homotopy between the two homotopies $(H_X)|_Y$ and 
	$\aof{l} \circ H_Y$
	in the category of spectra.  Hence there is a natural homotopy between 
	$\atrans{p_{X,Y}}\Psi(\alpha)$ and
	$\atransdec{p_{X,Y}} \Psi(\alpha)$.
\end{prop}

\begin{proof}
	By natural homotopy commutativity of the diagram given by  Lemma \ref{techone}
	and the various lifts established in
	Lemma \ref{samecomp}, it will suffice to see there is a natural homotopy 
	in $\aofp[Y]$ between $Res(\gamma_X)$ and $\omega_Y$.  However, this is a consequence
	of the construction of $\omega_Y$ from $Res(\gamma_X)$ and two other natural paths.
	
	The second claim follows from the first by employing 
	Def. \ref{makepairs} and the 
	compatibility of the homotopies established in the first claim.
\end{proof}

We now have a pair of  technical lemmas which will allow us to verify the natural homotopy commutativity of our key diagram.  First we have a consequence of the the Bousfield-Kan \cite{BK} models for homotopy limits and colimits.

\begin{lemma}
	\label{techtwo}
	Suppose $F:\catc \to \spectra$ and $\varphi:\catd \to \catc$ are functors, with $\catc$ and $\catd$ 	small categories.  Then the following diagram commutes up to a natural homotopy 
\def\objectstyle{\scriptstyle}
\diagramit{
{\holimdown{\catc} F} \ar[r]^{\varphi^*} \ar[d]_{\rho_\catc} & 
            {\holimdown{\catd} F \circ \varphi} \ar[d]^{\rho_\catd} \\
{\mapsp{\hocolimdown\catc |\catc/?|,\hocolimdown{\catc} F} } \ar[d]_{\epsilon}  &  
	   {\mapsp{\hocolimdown\catd |\catd/?|,\hocolimdown{\catd} F \circ \varphi}} \ar[dl]^{\delta} \\
 {\mapsp{\hocolimdown\catd |\catd/?|,\hocolimdown{\catc} F}},
}	
where $\epsilon$ denotes precomposition by 
\[\hocolim \varphi/?:\hocolimdown\catd |\catd/?| \to \hocolimdown\catc |\catc/?|
\] 
and $\delta$ denotes postcomposition by
\[\hocolimdown{\catd} F \circ \varphi \to \hocolimdown{\catc} F.
\]
\end{lemma}

\begin{lemma}
	\label{techone}
	The diagram 
\diagramit{
{\aofp[X]} \ar[d]_{\nu^{E_X}} \ar[rr]^{Res} 
   && {\aofp[Y]} \ar[d]^{\nu^{E_Y}} \\
{\holimdown{\tau \in \simp X} \aof{E_X^\tau}} \ar[rr] \ar[d]_{\rho^{E_X}} && 
            {\holimdown{\sigma \in \simp Y}\aof{E_Y^\sigma}} \ar[d]^{\rho^{E_Y}} \\
{\mapsp{X,\aof{E_X}} } \ar[dr]_{\epsilon}  &&  
	   {\mapsp{Y,\aof{E_Y}}} \ar[dl]^{\delta} \\
& {\mapsp{Y,\aof{E_X}}}
}		
	commutes up to a natural homotopy in the category of spectra.
\end{lemma}

\begin{proof}
	The homotopy commutativity of the top square is just an observation about the naturality of
	the Thomason homotopy limit problem map 
\[\nu_*^{E_X}:\aofp[X] \to \holimdown{\tau \in \simp X} \aof{E_X^\tau},
\] 
	as discussed by Williams \cite[page 12]{Wil}.

	The lower pentagon is in part an application of Lemma \ref{techtwo}, with $\catc=\simp{X}$, 
	$\catd=\simp{Y}$, $\varphi$ the map induced by the topological inclusion $Y \to X$ and
	$F$ the value of $\aof{?}$ at the pullback of $E_X$ over the inclusion of a given simplex, or
	explicitly $F(\sigma)=\aof{E_X^\sigma}$.  Notice if $\sigma=\varphi(\tau)$, we see 
	$E_X^\sigma\approx E_Y^\tau$ by construction, hence
\[
F \circ \varphi(\tau)=\aof{E_X^\sigma}= \aof{E_Y^\tau}.
\] 
	Thus, we have the natural homotopy commutativity of
\diagramit{
{\holimdown{\tau \in \simp X} \aof{E_X^\tau}} \ar[rr] \ar[d]_{\rho^{E_X}} && 
            {\holimdown{\sigma \in \simp Y}\aof{E_Y^\sigma}} \ar[d]^{\rho^{E_Y}} \\
{\mapsp{X,\abof[X]{E_X}} } \ar[dr]_{\epsilon}  &&  
	   {\mapsp{Y,\abof[Y]{E_Y}}} \ar[dl]^{\delta} \\
& {\mapsp{Y,\abof[X]{E_X}}}
}	

	Now the following diagram commutes, where the vertical maps are induced by
	$u^{E_X}$ (or $u^{E_Y}$)	
\diagramit{
{\mapsp{X,\abof[X]{E_X}} } \ar[r] \ar[d] & {\mapsp{Y,\abof[X]{E_X}}} 
     \ar[d] & {\mapsp{Y,\abof[Y]{E_Y}}} \ar[l] \ar[d] \\
{\mapsp{X,\aof{E_X}} } \ar[r] & {\mapsp{Y,\aof{E_X}}} & {\mapsp{Y,\aof{E_Y}}} \ar[l] 
}		
	which completes the proof.
\end{proof}

\end{document}